\documentclass[12pt]{amsart}

\usepackage{amssymb,amsmath,amsthm,amscd}
\usepackage{latexsym,amsfonts,epsfig,verbatim}
\usepackage{color,graphics, textcomp}

\input{epsf.tex}

\newcommand{\nc}{\newcommand}
\nc{\dmo}{\DeclareMathOperator}

\dmo{\mcg}{Mod}
\dmo{\pmcg}{PMod}

\def\mod{\mcg}
\def\mods{\mcg(S)}

\nc{\son}{S_{1,n}}
\nc{\modson}{\mcg(\son)}

\nc{\B}{\mathcal{A}}
\nc{\C}{\mathcal{C}}
\nc{\N}{\mathcal{N}}

\dmo{\comm}{Comm}
\dmo{\aut}{Aut}

\dmo{\gin}{i}

\newtheorem{thm}{Theorem}

\newtheorem{main}[thm]{Main Theorem}
\newtheorem{prop}[thm]{Proposition}

\newtheorem*{qn}{Question}
\newtheorem{lem}[thm]{Lemma}

\theoremstyle{definition}

\theoremstyle{remark}

\nc{\Z}{\mathbb{Z}}
\nc{\R}{\mathbb{R}}
\nc{\pd}{\mathcal{P}}
\nc{\G}{\mathcal{G}}
\nc{\T}{\mathcal{T}}

\nc{\al}{\alpha}
\nc{\be}{\beta}
\nc{\ga}{\gamma}

\nc{\p}[1]{\paragraph{{\bf #1}}}
\nc{\set}[1]{\{#1\}}
\nc{\bpf}{\begin{proof}}
\nc{\epf}{\end{proof}}
\nc{\margin}[1]{\marginpar{\scriptsize #1}}

\nc{\pics}[3]{\epsfysize=#3 cm \begin{figure}[htb]
\center{\leavevmode \epsfbox{#1.eps}} \caption{#2}  \label{#1pic}
\end{figure}}

\begin{document}

\title[Curve complexes and subgroups of mapping class groups]{Curve complexes and finite index subgroups of mapping class groups}

\author{Jason Behrstock}
\address{Department of Mathematics \\ Barnard College,
Columbia University \\ New York, NY 10027}
\email{jason@math.columbia.edu}

\author{Dan Margalit}

\address{Department of Mathematics\\ University of Utah\\ 155 S 1440
East \\ Salt Lake City, UT 84112-0090}

\thanks{The second author is partially supported by an NSF
postdoctoral fellowship.}

\email{margalit@math.utah.edu}



\maketitle

\begin{center}\today\end{center}

\begin{abstract} Let $\mods$ be the extended mapping class group of a surface $S$.  For $S$ the twice-punctured torus, we show that there exists an isomorphism of finite index subgroups of $\mods$ which is not the restriction of an inner automorphism.  For $S$ a torus with at least three punctures, we show that every injection of a finite index subgroup of $\mods$ into $\mods$ is the restriction of an inner automorphism; this completes a program begun by Irmak.  For all of the above surfaces, we establish the co-Hopf property for finite index subgroups of $\mods$.\end{abstract}

\section{Introduction} \label{intro}

Let $S=S_{g,n}$ be a \emph{surface}, by which we always mean a connected, orientable surface of genus $g$ with $n$ punctures.  The \emph{extended mapping class group} of $S$ is:
\[ \mods = \pi_0 (\textrm{Homeo}^{\pm}(S)) \]
In this note we complete the answer to the following question:

\begin{qn}Given a surface $S$, is it true that every
injection of a finite index subgroup of $\mods$ into $\mods$ is the restriction of an inner automorphism?\end{qn}

The study of maps of finite index subgroups of $\mods$ into $\mods$
was begun by Ivanov, although this specific question was first
investigated by Irmak, who showed that the answer is no for $S_{2,0}$
(see below) and yes for all other $S_{g,n}$ with $g \geq 2$ \cite{ei}
\cite{ei2} \cite{ei3}.  Bell--Margalit proved that the answer is yes
for $S_{0,n}$ with $n \geq 5$, giving applications to the study of
Artin groups \cite{bm}.

Korkmaz showed that the answer is no for $S_{0,4}$, $S_{1,1}$, and $S_{1,0}$, as the mapping class groups for these surfaces contain free subgroups with finite index \cite{mkcof}. Also, the answer is no for $S_{0,2}$ and $S_{0,3}$ and yes for $S_{0,0}$ and $S_{0,1}$; this is not hard to check since the mapping class groups for these surfaces are finite.

In this paper, we treat all remaining cases.  For $S_{1,2}$, we show that the answer is no:

\begin{main} \label{main1} There exists an isomorphism between finite index subgroups of $\mod(S_{1,2})$ which is not the restriction of an inner automorphism of $\mod(S_{1,2})$. \end{main}

We show that the answer is yes in the other cases:

\begin{main} \label{main2} Let $S$ be $S_{1,n}$ with $n \geq 3$.  Any injection of a finite index subgroup of $\mods$ into $\mods$ is the restriction of an inner automorphism. \end{main}

We thus have the following theorem, which combines our results with those of Korkmaz, Irmak, and Bell--Margalit.

\begin{thm} \label{injthm} If $S$ is a surface which is not $S_{0,2}$, $S_{0,3}$, $S_{0,4}$, $S_{1,0}$, $S_{1,1}$, $S_{1,2}$, or $S_{2,0}$, then every injection of a finite index subgroup of $\mods$ into $\mods$ is the restriction of an inner automorphism of $\mods$. If $S$ is one of these exceptional surfaces, then there is an isomorphism of finite index subgroups of $\mods$ which is not the restriction of an inner automorphism. \end{thm}

In the case of $S_{2,0}$, if $\Gamma$ is a finite index subgroup of
$\mcg(S_{2,0})$ and $\rho:\Gamma \to \mcg(S_{2,0})$ is an injection,
then $\rho$ is given by $\rho(g) = fgf^{-1} \iota^{\sigma(g)}$, where
$f \in \mcg(S_{2,0})$, $\iota$ is the hyperelliptic involution, and $\sigma: \Gamma
\to \Z_2$ is a homomorphism; note $\langle \iota \rangle =
Z(\mcg(S_{1,2})) \cong \Z_2$.  It is not hard to construct examples of injections where the associated homomorphism $\sigma$ is nontrivial.  If $\Gamma'$ is the kernel of $\sigma$,
then $\rho|_{\Gamma'}$ is the restriction of an inner automorphism and the index of $\Gamma'$ in $\Gamma$ is at most 2.

\p{Complex of curves and Ivanov's theorem.} Let $\C(S)$ denote the \emph{complex of curves} for $S$, which is the abstract simplicial flag complex with a vertex for each isotopy class of simple closed curves in $S$ and an edge between vertices with disjoint representatives; this complex was defined by Harvey \cite{wjh}.  In his seminal work, Ivanov proved that every isomorphism between finite index subgroups of $\mods$ is the restriction of an inner automorphism; the main step was to show that every automorphism of $\C(S)$ is induced by an element of $\mods$ \cite{nvi}.

In light of Ivanov's theorem, Theorem~\ref{injthm} can be thought of as saying that every injection of a finite index subgroup of $\mods$ has finite index image in $\mods$.  However, showing that the image of such an injection has finite index does not seem to be easier than showing directly that the injection is the restriction of an inner automorphism.

\p{Superinjective maps.} To attack the particular question at hand, Irmak introduced the notion of a \emph{superinjective} map of $\C(S)$, which is a simplicial map of $\C(S)$ to itself preserving disjointness and nondisjointness of the isotopy classes of curves corresponding to the vertices of $\C(S)$ (note any simplicial map of $\C(S)$ preserves disjointness).

It is clear that elements of $\mods$ give rise to superinjective maps of $\C(S)$.  Following Irmak, Main Theorem~\ref{main2} is proven by showing that all superinjective maps arise in this way:

\begin{thm}\label{sithm} Let $S$ be $S_{1,n}$ with $n \geq 3$.  Every superinjective map of $\C(S)$ is induced by an element of $\mods$. \end{thm}

For this theorem, the main difficulty lies in
distinguishing nonseparating curves from curves bounding
twice-punctured disks.  To overcome this obstacle, we introduce an essential new tool, the \emph{adjacency graph} of a pants
decomposition of a surface (see Section~\ref{g1sec}).
This notion can also be used to provide a new approach
to Theorem~\ref{sithm} for the cases not covered in this paper.

After the completion of this work, we learned that the idea of the
adjacency graph was independently (and simultaneously) discovered by Shackleton.  He used it to show (in most cases) that any simplicial embedding of a curve complex into a curve complex of equal or lesser dimension is necessarily an automorphism \cite{kjs}.  This reproves Theorem~\ref{injthm} in many cases, and also has the consequence that,
for these mapping class groups, there are no injections of a finite index subgroup into a mapping class group of equal or lesser complexity.

The deduction of Main Theorem~\ref{main2} from Theorem~\ref{sithm} is now a standard argument for which we refer the reader to Irmak's paper \cite{ei}.  The idea is that an injection of a finite index subgroup of $\mods$ into $\mods$ must take powers of \emph{Dehn twists} to powers of Dehn twists.  These are exactly the elements of $\mods$ which are each supported on the regular neighborhood of a simple closed curve in $S$ (see e.g. \cite{jb}).  Powers of Dehn twists commute if and only if the corresponding curves are disjoint, so the injection gives a superinjective map of $\C(S)$.

\p{Conventions.} When there is no confusion, we will use \emph{curve} to mean ``isotopy class of an essential simple closed curves,'' and we will blur the distinction between a curve and its isotopy class. When we say two curves $a$ and $b$ are \emph{disjoint}, we mean that their geometric intersection number $\gin(a,b)$ is 0. An \emph{arc} is the isotopy
class of the image of a proper, essential embedding of $\R$ in a surface;
here, \emph{essential} means not homotopic to a puncture. 

\p{Acknowledgments.} The first author dedicates this work to his
daughter, Kaia Behrstock, who was born at the onset of this
paper. Both authors would like to extend their thanks to the Columbia University Math Department, where the second author was visiting at the origin of this project.  We thank Bob Bell, Mladen Bestvina, Richard Kent, Chris Leininger, Ken Shackleton, and Steven Spallone for being very generous with their time and energy.  We appreciate Elmas Irmak's comments on an earlier draft.  We are also grateful to Kashi Behrstock, Benson Farb, and Kevin Wortman for their comments and encouragement.


\section{The co-Hopf property}

A group $G$ is \emph{co-Hopfian} if every injective endomorphism of $G$ is an automorphism of $G$, and it is \emph{Hopfian} if every surjective endomorphism is an automorphism.  In general, the former property seems to be more rare and harder to prove than the latter.

Grossman proved that mapping class groups are
residually finite \cite{ekg} (see also \cite{ni}).  Since finitely generated residually finite groups are Hopfian, we have this property for $\mods$.

It follows formally from Theorem~\ref{injthm} that $\mods$ is
co-Hopfian for most surfaces $S$.  This was first proven by
Ivanov--McCarthy \cite{im}.  What is more, it is straightforward to deduce from Theorem~\ref{injthm} the following:

\begin{thm}\label{cohopfthm}If $S$ is not $S_{0,4}$, $S_{1,0}$, or $S_{1,1}$, then every finite index subgroup of $\mods$ is co-Hopfian.\end{thm}

Note that the theorem is not true for $S$ either $S_{0,4}$, $S_{1,0}$,
or $S_{1,1}$, by the fact that $\mods$ contains a free group with
finite index in these cases.  Since they are not treated by
Theorem~\ref{injthm}, the surfaces $S_{2,0}$ and $S_{1,2}$ require
special attention.  The case of $S_{2,0}$ is handled by the statement
immediately following Theorem~\ref{injthm} and an elementary group
theory argument.  We relegate the argument for $S_{1,2}$ to Section~\ref{sotcohopf}.

\p{Related results.}  Farb--Ivanov showed that the \emph{Torelli group}, the subgroup of $\mods$ acting trivially on the homology of $S$, is co-Hopfian \cite{fi}.  Brendle--Margalit proved that the so-called \emph{Johnson kernel} and all of its finite index subgroups are co-Hopfian \cite{brm};  this result and the previous were proven for closed surfaces of genus at least 4.
Bell--Margalit established the co-Hopf property for the \emph{braid group on n strands} $B_n$ modulo its center when $n \geq 4$ \cite{bm1}; this is essentially the genus 0 version of the Ivanov--McCarthy theorem \cite{bm}. Recently, Farb--Handel showed that $\textrm{Out}(F_n)$ and all of its finite index subgroups are co-Hopfian for $n \geq 4$ \cite{fh}.


\section{Abstract commensurators}

To better understand Main Theorem~\ref{main1}, it will be helpful to recast Ivanov's theorem by putting a group structure on the set of isomorphisms of finite index subgroups of $\mods$.

The \emph{abstract commensurator} $\comm(G)$ of a group $G$ is the group of equivalence classes of isomorphisms of finite index subgroups of $G$.  Two isomorphisms are said to be equivalent if they agree on a finite index subgroup of $G$.  The composition of two isomorphisms $\psi: \Gamma \to \Lambda$ and $\psi':\Gamma' \to \Lambda'$ is a map defined on the finite index subgroup $\psi^{-1}(\Lambda \cap \Gamma')$.  A simple example is $\comm(\mathbb{Z}^n) \cong \textrm{GL}_n(\mathbb{Q})$.

In this language, Ivanov's result is:

\begin{thm}\label{ivanovthm}Let $S$ be any surface other than
  $S_{1,n}$ with $n \leq 2$ or $S_{0,n}$ with $n \leq 4$.  The natural
  homomorphism \[ \mods \to \comm(\mods) \] is surjective.  Moreover,
  the kernel is $Z(\mods) \cong \Z_2$ for $S=S_{2,0}$ and is trivial otherwise.\end{thm}

This theorem was proven by Ivanov for surfaces of genus at least 2 \cite{nvi}.  Korkmaz then proved it for surfaces of genus 0 and 1
\cite{mk}.  Theorem~\ref{ivanovthm} can also be deduced from work of Luo, who gave a new proof that the natural map $\mods \to \aut(\C(S))$ is surjective \cite{fl}.

We remark that Theorem~\ref{ivanovthm} is a corollary of
Theorem~\ref{injthm}.  To see why Theorem~\ref{injthm} is a priori
much more difficult, note that, by the same logic used to deduce Main
Theorem~\ref{main2} from Theorem~\ref{sithm} (see
Section~\ref{intro}), an element of $\comm(\mods)$ gives rise to an
\emph{automorphism} of $\C(S)$.  For the non-exceptional surfaces, the
topological type of a curve is determined by the homotopy type of its
link in $\C(S)$ (this is a theorem of Harer \cite{jlh}).  It follows
immediately that an element of $\aut(\C(S))$ preserves the topological
type of a curve.  For a superinjective map of $\C(S)$, this is a
difficult step.  Still, our argument follows the general outline of
Ivanov's proof.

\p{Related results.}  By work of Farb--Ivanov and Brendle--Margalit, the abstract commensurators of the Torelli group and the Johnson kernel are both isomorphic to $\mods$ for closed surfaces of genus at least 4 \cite{fi} \cite{brm}.
Charney--Crisp used Theorem~\ref{ivanovthm} to show that the abstract
commensurators of several affine and finite type Artin groups, modulo
their centers, are isomorphic to $\mcg(S_{0,n})$ \cite{cc}.
Leininger--Margalit proved that $\comm(B_n) \cong \mcg(S_{0,n+1}) \ltimes (\mathbb{Q}^\times \ltimes
\mathbb{Q}^\infty)$ when $n \geq 4$ \cite{lm}, and Farb--Handel showed
$\comm(\textrm{Out}(F_n)) \cong \textrm{Out}(F_n)$ for $n \geq 4$ \cite{fh}.


\section{Twice-punctured torus} \label{sotsec}

Luo showed that $S_{1,2}$ has the exceptional property that there are
automorphisms of $\C(S_{1,2})$ which are not induced by
$\mod(S_{1,2})$ \cite{fl}.  We will translate this fact into the group
theoretic statement of Main Theorem~\ref{main1}.  We then establish
Theorem~\ref{cohopfthm} for $S_{1,2}$, and also describe all
injections of finite index subgroups of $\mod(S_{1,2})$ into $\mod(S_{1,2})$.

\subsection{Nongeometric curve complex automorphisms.} We give the
argument of Luo.  The first step is to construct an isomorphism between $\C(S_{0,5})$ and $\C(S_{1,2})$.  Let $\iota$ be the hyperelliptic involution of $S_{1,2}$.  The quotient $S_{1,2}/\langle \iota \rangle$ is a sphere with 1 puncture and 4 cone points of order 2. Identify $S_{0,5}$ with the complement of a regular neighborhood of these cone points.  Given any (simple closed) curve in $S_{0,5}$, we can lift it via $\iota$ to a curve in $S_{1,2}$.  Moreover, this map preserves disjointness. Also, this map is surjective: Birman and Viro showed that $\iota$ fixes every curve in $S_{1,2}$ (being central, it commutes with each Dehn twist), and so there is an inverse map \cite{jb} \cite{ojv}.

In $S_{0,5}$, all curves are topologically equivalent, since they each have three punctures on one side and two on the other.  Thus, $\aut(\C(S_{0,5}))$ acts transitively on the vertices of $\C(S_{0,5})$.  Since $\C(S_{1,2})$ is isomorphic to $\C(S_{0,5})$, it follows that there are elements of  $\aut(\C(S_{1,2}))$ which interchange separating and nonseparating curves.  Such automorphisms clearly cannot be induced by $\mod(S_{1,2})$.

To see this more concretely, note that the curves in $S_{0,5}$ corresponding to separating curves in $S_{1,2}$ are exactly the curves which contain the \emph{special puncture} (the one coming from the punctures of $S_{1,2}$) on their twice-punctured sides.

\subsection{Nongeometric commensurators.} We will show that the elements of $\aut(\C(S_{1,2}))$ which are not induced by $\mod(S_{1,2})$ give rise to elements of $\comm(\mod(S_{1,2}))$ which are not induced by $\mod(S_{1,2})$.

Let $\pmcg(S_{1,2})$ denote the subgroup of $\mod(S_{1,2})$ consisting of elements which fix each puncture.  By work of Birman--Hilden, $\iota$ induces an isomorphism $\iota_\star$ of $\pmcg(S_{1,2})$ with the subgroup of $\mcg(S_{0,5})$ consisting of elements which fix the special puncture \cite{bh72}.  This index 5 subgroup of $\mcg(S_{0,5})$ corresponds exactly to the subgroup of $\aut(\C(S_{0,5}))$ consisting of elements which lift to automorphisms of $\C(S_{1,2})$ induced by $\mcg(S_{1,2})$.

Combining Theorem~\ref{ivanovthm} with the fact that $\mod(S_{1,2})$ and $\mcg(S_{0,5})$ have isomorphic finite index subgroups, it follows that:

\begin{prop} \label{comms}
$\comm(\mod(S_{1,2})) \cong \comm(\mcg(S_{0,5})) \cong \mcg(S_{0,5})$.
\end{prop}

Now, consider the composition
\[ \mod(S_{1,2}) \to \comm(\mod(S_{1,2})) \to \comm(\mcg(S_{0,5})) \to \mcg(S_{0,5}) \]
where $g \in \mcg(S_{1,2})$ maps to conjugation by $g$, which then maps to conjugation by $\iota_\star(g)$, and finally to the element $\iota_\star(g)$ of $\mcg(S_{0,5})$. We have already noted that $\iota_\star(\mod(S_{1,2}))$ has index 5 in $\mcg(S_{0,5})$.  Since the kernel of the first map is $Z(\mod(S_{1,2})) = \langle \iota \rangle \cong \Z_2$ and the last two maps are the restrictions of the isomorphisms of Proposition~\ref{comms}, this yields the following, which proves Main Theorem~\ref{main1}.

\begin{prop}$[\comm(\mod(S_{1,2})):\mod(S_{1,2})/\Z_2]=5.$ \end{prop}

We now explain why the elements of $\comm(\mcg(S_{1,2}))$ which do not come from $\mcg(S_{1,2})$ are not geometric in any sense.

If $\pmcg(S_{0,5})$ is defined similarly to $\pmcg(S_{1,2})$, then $\pmcg(S_{0,5})$ is normal in $\mcg(S_{0,5})$ and is isomorphic (via $\iota_\star$) to a finite index subgroup $\Gamma$ of $\mod(S_{1,2})$.  Now, we can write the isomorphism
\[ \Upsilon : \mcg(S_{0,5}) \to \comm(\mcg(S_{1,2})) \]
as:
\[ g \mapsto [\Gamma \to \iota_\star^{-1}(g \iota_\star(\Gamma) g^{-1}) ] \]
Denote by $T_c$ the Dehn twist about a curve $c$, and let $g$ be an arbitrary element of $\mods$.  From the general formula
\begin{equation}\label{conjtwist} gT_c^kg^{-1}=T_{g(c)}^{\pm k} \end{equation}
it follows that any representative of $\Upsilon(g)$ takes powers of $T_c$ to powers of $T_{g_\star(c)}$, where $g_\star$ is the induced element of $\aut(\C^0(S_{1,2}))$ ($\C^0(S_{1,2})$ is the 0-skeleton).  Since any element of $\mcg(S_{0,5}) - \iota_\star(\mcg(S_{1,2}))$ gives an element of $\aut(\C^0(S_{1,2}))$ interchanging separating and nonseparating
curves, it follows that the image under $\Upsilon$ of these elements
takes powers of Dehn twists about separating curves to powers of Dehn
twists about nonseparating curves (and vice versa).  By contrast,
equation~(\ref{conjtwist}) shows that an inner automorphism preserves
the topological types of curves corresponding to powers of Dehn twists.

\p{Remark.}  The elements of $\comm(\mcg(S_{1,2}))$ which do not arise from inner automorphisms of $\mcg(S_{1,2})$ actually do not arise from any elements of $\aut(\mcg(S_{1,2}))$.  This follows from the fact that, in $\mcg(S_{1,2})$, conjugate Dehn twists can commute only if they are twists about nonseparating curves.  Thus,  $\aut(\mcg(S_{1,2})) \cong \mcg(S_{1,2})/\langle \iota \rangle$.

\subsection{Proof of Theorem~\ref{cohopfthm}}\label{sotcohopf} Despite the fact that there are nongeometric elements of $\comm(\mcg(S_{1,2}))$, we can still establish Theorem~\ref{cohopfthm} in the case of $S_{1,2}$.  The first step is:

\begin{prop}\label{direct} $\mcg(S_{1,2}) \cong \pmcg(S_{1,2}) \times \mathbb{Z}_2$.\end{prop}

This proposition follows from the fact that the short exact sequence
\[ 1 \to \langle \iota \rangle \to \mcg(S_{1,2}) \to \mcg(S_{1,2})/\langle \iota \rangle \to 1 \]
has a splitting $\mcg(S_{1,2})/\langle \iota \rangle \to \mcg(S_{1,2})$ defined by sending the coset $g \langle \iota \rangle$ to its unique representative in $\pmcg(S_{1,2})$.

Since $\pmcg(S_{1,2})$ is isomorphic to a finite index subgroup of
$\mcg(S_{0,5})$, every finite index subgroup of $\pmcg(S_{1,2})$ is
co-Hopfian, as Theorem~\ref{cohopfthm} follows from
Theorem~\ref{injthm} in the case of $S_{0,5}$.  Using
Proposition~\ref{direct} and the fact that finite index subgroups of
$\mcg(S_{0,5})$ have trivial center, an elementary group theory
argument shows that every finite index subgroup of $\mcg(S_{1,2})$ is
co-Hopfian.

\subsection{Injections of finite index subgroups}

If $K$ is any finite index subgroup of $\pmcg(S_{1,2}) \times \Z_2
\cong \mcg(S_{1,2})$, then $K$ is isomorphic to $G \times H$, where
$G$ is a finite index subgroup of $\pmcg(S_{1,2})$ and $H < \Z_2$.  We
fix inclusions:
\[ G \hookrightarrow \pmcg(S_{1,2}) \hookrightarrow \mcg(S_{0,5}) \]
By an argument similar to that needed for Section~\ref{sotcohopf}, we can describe all injections of finite index subgroups of $\mcg(S_{1,2})$ into $\mcg(S_{1,2})$:

\begin{prop}With the above notation, let $\rho$ be an injection $K \to \mcg(S_{1,2})$, thought of as $G \times H \to \pmcg(S_{1,2}) \times \Z_2$.  Then $\rho$ is given by
\[ (g,x) \mapsto (\Omega(g),\sigma(g) \cdot x) \]
where $\Omega \in \aut(\mcg(S_{0,5})) \cong \mcg(S_{0,5})$ and $\sigma: G \to \Z_2$ is a homomorphism.\end{prop}

When, in the conclusion of the proposition, $\Omega$ is trivial, we
remark that $\rho$ is a \emph{transvection} in a sense slightly generalizing that
of Charney--Crisp \cite{cc} (see also \cite{bm} \cite{lm}).


\section{Proof of Theorem~\ref{sithm}} \label{g1sec}

In this section, we fix a superinjective map
$\phi:\C(\son) \to \C(\son)$, where $n \geq 3$.
We will prove Theorem~\ref{sithm} in this case; that is, we will show
that $\phi$ is induced by an element of $\modson$.

We start by showing that $\phi$ preserves certain topological
properties of (and relationships between) curves.

\p{Pants decompositions.} It is not hard to see that any
superinjective map is injective (for any two curves, consider a curve
which is disjoint from one but not the other).  It follows that
$\phi$ preserves the set of  \emph{pants decompositions} of $\son$,
i.e. the maximal simplices of $\C(\son)$.

\p{Adjacency graphs.} Given a pants decomposition $\pd$ of $\son$, we
say that two curves of $\pd$ are \emph{adjacent with respect to $\pd$}
if there exists a component of $\son-\pd$ containing both curves in
its closure; this definition is essential in the work of Irmak
\cite{ei} \cite{ei2}.  We assign an \emph{adjacency graph} $\G(\pd)$
to $\pd$ consisting of a vertex for each curve in $\pd$ and edges
corresponding to  adjacency.

The following result provides our main
tool, and the proof we give works for all surfaces.

\begin{lem}\label{graphspreserved}
For any pants decomposition $\pd$, $\phi$ induces an isomorphism of
the graphs $\G(\pd)$ and $\G(\phi(\pd))$.
\end{lem}

\bpf It suffices to show $\phi$ preserves adjacency and nonadjacency
with respect to $\pd$.

Two curves are adjacent with respect to $\pd$ if and only if there
exists a curve which intersects both and is disjoint from the other
curves of $\pd$.  Since superinjective maps preserve disjointness and
nondisjointness, this characterization implies that $\phi$ preserves
adjacency.

Say that the curves of $\pd$ are $a_1,b_1,c_1, \dots, c_{n-2}$.  If
$a_1$ and $b_1$ are not adjacent with respect to $\pd$, then there
exist curves $a_2$ and $b_2$ so that $\set{a_i,b_j,c_1, \cdots,
c_{n-2}}$ is a pants decomposition for any $i,j \in \set{1,2}$.
Since curves in a pants decomposition are disjoint, it follows that
$a_1$, $b_1$, $a_2$, and $b_2$ (in that order) form a square in
$\C(\son-\cup c_i)$.

Now, if $a_1$ and $b_1$ are in fact adjacent with respect to $\pd$,
then they lie on a common connected component of $\son-\cup c_i$,
which must be homeomorphic to $S_{0,5}$ or $S_{1,2}$.  As there are
no squares in $\C(S_{0,5})$ or $\C(S_{1,2})$, it follows that there
are no curves $a_2$ and $b_2$ with the property that
$\set{a_i,b_j,c_1, \cdots, c_{n-2}}$ is a pants decomposition for any
choice of $i,j \in \set{1,2}$.

Since $\phi$ must preserve the existence of the given squares in
$\C(\son-\cup c_i)$, it follows that superinjectivity must preserve
the property of nonadjacency.
\epf

To simplify the arguments, we will make the assumption of $n \geq 4$ in Lemmas~\ref{linearorcyclic} through~\ref{dualitypreserved}.  After Lemma~\ref{dualitypreserved}, we explain how the argument should be adjusted for the case of $n=3$.

\p{Linear and cyclic pants decompositions.} Our next goal is to show
that $\phi$ preserves the topological types of curves.  The first step
is to give
a classification of pants decompositions of $\son$. Below, a
\emph{linear} pants decomposition is one which is topologically
equivalent to the one on the left hand side of Figure~\ref{lincycpic}
and a \emph{cyclic} pants decomposition is one which is equivalent to
the one on the right (the names come from the shapes of the
corresponding adjacency graphs).  Note that since $n \geq 3$, the
adjacency graphs are distinct.

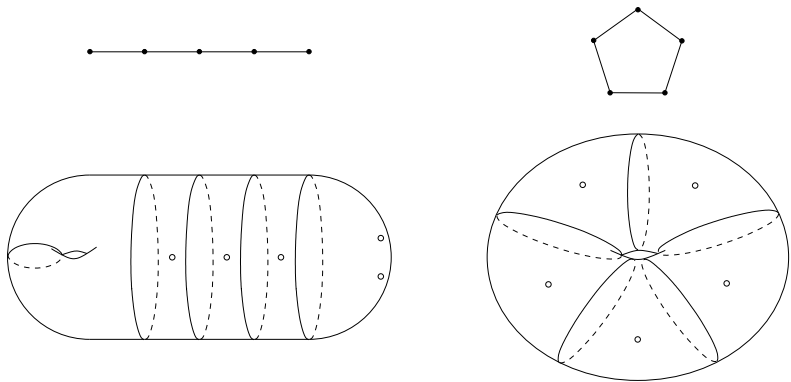
\pics{lincyc}{Linear (left) and cyclic (right) pants decompositions
for $S_{1,5}$.}{5}

We shall call the closures of the components of $\son-\pd$ the
\emph{pairs of pants} of $\pd$.

\begin{lem}\label{linearorcyclic}
Let $ n\geq 4$.  If $\pd$ is a pants decomposition of $\son$ with the property that $\G(\pd)$ does not contain a triangle, then $\pd$ is either a linear pants decomposition or a cyclic pants decomposition.
\end{lem}

\bpf None of the pairs of pants of $\pd$ can have 3 distinct curves of
$\pd$ for their boundary, since that would form a triangle in the
adjacency graph.  Thus, any
pair of pants of $\pd$ is of one of the following types:
\begin{enumerate}
    \item torus with one boundary
    \item punctured annulus
    \item twice-punctured disk
\end{enumerate}

We imagine constructing $\son$ out of these pieces.  The only
surfaces without boundary which can be built from pairs of pants of
types 1 and 3 are $S_{2,0}$, $S_{0,4}$, and $S_{1,2}$.  Thus, we may
conclude that $\pd$ has at least one pair of pants of type 2.
Moreover, all of the type 2 pairs of pants of $\pd$ lie in a
connected chain, since each pair of pants of type 1 or 3 is connected
to only one other pair of pants.

If we add a type 1 pair of pants to each end of this chain of type 2
pairs of pants, then the result is $S_{2,m}$ for some $m$.  If we cap
off the chain by a two pairs of pants of type 3, the result is
$S_{0,m}$ for some $m$.  We can attach one type 1 and one type 3 pair
of pants to the ends of the chain---this yields a linear pants
decomposition of $\son$.  If we fail to attach a pair of pants to
either end of this chain, we get a surface with boundary, so the only
other possibility is to attach the ends of the chain together---this
gives a cyclic pants decomposition of $\son$.
\epf

A \emph{side} of a separating curve $z$ is a connected component of
$\son-z$.  In the following lemma, when we say $\phi$ \emph{preserves
sides} of a separating curve $z$, we mean that if $a$ and $b$ are curves on the same side of $z$, then $\phi(a)$ and $\phi(b)$ are on the same side of
$\phi(z)$.

\begin{lem}\label{Npreserved}Let $n \geq 4$. $\phi$ preserves the topological types of curves, and $\phi$ preserves sides of separating curves.
\end{lem}

\bpf
First, if $a$ is a nonseparating curve in $\son$, then $a$ fits into
a cyclic pants decomposition $\pd$ of $\son$.  By
Lemmas~\ref{graphspreserved} and~\ref{linearorcyclic}, $\phi(\pd)$ is
cyclic.  Since all curves in a cyclic pants decomposition are
nonseparating, $\phi(a)$ must be nonseparating.

Now if $a$ is any curve in $\son$, then $a$ fits into a linear pants
decomposition $\pd$ of $\son$.  Again, by Lemmas~\ref{graphspreserved}
and~\ref{linearorcyclic}, $\phi(\pd)$ is linear, and $\phi$ induces an
isomorphism of $\G(\pd)$ with $\G(\phi(\pd))$.  We have already shown
that $\phi$ preserves nonseparating curves, so $\phi$ must take the
unique vertex of $\G(\pd)$ representing a nonseparating curve to the
unique vertex of $\G(\phi(\pd))$ representing a nonseparating curve.
Since the rest of the graph isomorphism is determined, it follows
that $\phi$ preserves the topological type of every curve in $\pd$,
in particular $a$.

It follows from the above argument that $\phi$ preserves sides of
separating curves.
\epf

\p{Duality.} By a \emph{$k$-curve}, we will mean a separating curve
in $\son$ with exactly $k$ punctures on its genus 0 side.  Denote
by $\N$ the set of all nonseparating curves and 2-curves in $\son$.

We call a pair of curves $a$ and $b$ of $\N$ \emph{dual} if they are
nonseparating and $\gin(a,b)=1$ or if they are 2-curves and
$\gin(a,b)=2$.

To show that duality is preserved by $\phi$, we will require the
following vocabulary: let $\N'$ denote the union of the 3-curves and
$n$-curves in $\son$.  A \emph{small} side of an element $z$ of $\N'$
is a side of $z$ homeomorphic to $S_{1,1}$ (in the case $z$ is an
$n$-curve) or $S_{0,4}$ (in the case $z$ is a 3-curve).  Note that when $n=3$, any element of $\N'$ has two small sides.

\begin{lem}\label{dualitypreserved}For $n \geq 4$, $\phi$ preserves duality.
\end{lem}

\bpf
Two curves $a$ and $b$ are dual if and only if there exists an
element $z$ of $\N'$ and curves $x$ and $y$ in $\son$ so that: $a$
and $b$ lie on the same small side of $z$, $x$ intersects $z$ and $a$ but
not $b$, and $y$ intersects $z$ and $b$ but not $a$.  One direction
is proven by construction.  The other direction is proven for
nonseparating curves by Ivanov and for 2-curves by Bell--Margalit
\cite{nvi} \cite{bm} (see also \cite{fl} \cite{mk}).

Now, note that the argument of Lemma~\ref{Npreserved} shows that $\phi$ preserves
small sides of curves.  Thus, $\phi$ preserves all of the properties
used in this characterization of duality, and so the lemma follows.
\epf

\p{The case $\boldsymbol{n=3}$.}  The subtle point when $n=3$ is that there are two topological types of pants decompositions whose adjacency graphs are triangles: those with either 2 or 3 nonseparating curves.  Thus, our argument in the proof of Lemma~\ref{Npreserved} that $\phi$ preserves nonseparating curves does not work.  However, since every pants decomposition of $S_{1,3}$ with a linear adjacency graph is linear in the sense of Lemma~\ref{linearorcyclic}, the argument of Lemma~\ref{Npreserved} does still show that $\phi$ preserves the sets $\N$ and $\N'$ (and sides of elements of $\N'$).  In that case, the argument of Lemma~\ref{dualitypreserved} shows that $\phi$ preserves duality in the case of $n=3$ (even though we don't know topological types are preserved!).  Now, we have the following distinction between nonseparating curves and 2-curves: given a nonseparating curve $c$, we can find three disjoint nonseparating curves which are all dual to $c$.  However, a 2-curve (on any surface) can only be dual to at most 2 disjoint curves.  Thus $\phi$ preserves nonseparating curves, and it then follows as in the proof of Lemma~\ref{Npreserved} that $\phi$ preserves 2-curves.

We remark that the above argument works in general for any $n \geq 3$.  However, it is more complicated than necessary when $n > 3$, and so we choose to single this case out.

We now return to the assumption that $n \geq 3$.  Below, when we refer to Lemmas~\ref{Npreserved} and~\ref{dualitypreserved}, we will mean the inclusive statements for all $n \geq 3$.

\bigskip

\p{Arc complex and ideal triangulations.} The \emph{arc complex} $\B(S)$ is the abstract simplicial flag complex with a vertex for each arc in the surface $S$ and edges corresponding to disjointness.

We define a \emph{triangle} to be a triple of pairwise dual 2-curves in $\son$ which lie on the genus 0 side of a 3-curve.  A \emph{2-curve ideal triangulation} is a collection of triangles corresponding to a maximal simplex in $\B(\son)$ (all maximal simplices of $\B(\son)$ have the same dimension).  It is not hard to check that there exist 2-curve ideal triangulations of $\son$ when $n \geq 3$.

We may think of a 2-curve ideal triangulation as an ideal triangulation in the usual sense by replacing each 2-curve with the unique arc connecting the punctures on the genus 0 side of the 2-curve.  Note that the arcs corresponding to two 2-curves are disjoint if and only if the 2-curves are disjoint or dual.

It follows immediately from Lemma~\ref{Npreserved},
Lemma~\ref{dualitypreserved}, the injectivity of $\phi$, and the definition of superinjectivity that we have:

\begin{lem}\label{triangulation} $\phi$ preserves triangles and 2-curve ideal triangulations.\end{lem}

If $\T$ is a 2-curve ideal triangulation of $\son$, then by Lemma~\ref{triangulation}, $\phi(\T)$ is a 2-curve ideal triangulation of $\son$ which is abstractly isomorphic to $\T$.  It follows that there is an element $f$ of $\mcg(\son)$ so that $f(\T)=\phi(\T)$.
To show that $f$ agrees with $\phi$ on the rest of $\C(\son)$, we will make use of the rigid combinatorial structure of $\B(S)$.

\p{Chain connectedness of $\boldsymbol{\B(S)}$.} It is not hard to see that any codimension 1 simplex of $\B(S)$ is contained in at most two simplices of maximal dimension. Also, there is a well-known theorem that $\B(S)$ is \emph{chain connected}, which means that any two maximal simplices are connected by a sequence of maximal simplices, where consecutive simplices in the sequence share a codimension 1 face.  There are various proofs of this theorem in the literature, by Mumford (see Harer's paper \cite{jlh}), Penner \cite{rcp}, Bowditch--Epstein \cite{be}, Mosher \cite{lm}, and Hatcher \cite{ah}.  Several of these papers credit the idea of this theorem to W. Thurston (see \cite{flp}).

From these facts about $\B(S)$, it follows that:

\begin{lem}\label{bs}Any injective simplicial map $\B(S) \to \B(S)$ is determined by its action on a single maximal simplex.\end{lem}

We already have that $\phi$ agrees with some element $f$
of $\mcg(\son)$ on a maximal simplex of $\B(\son)$, namely $\T$.  If
we show that $\phi$ extends to an injective simplicial map on all of
$\B(\son)$, it will follow from Lemma~\ref{bs} that $\phi$ agrees
with $f$ on all of $\B(\son)$.  From there, it will be easy to show that $\phi$ agrees with $f$ on $\C(\son)$.


\p{Induced map on $\boldsymbol{\B(\son)}$.} There is a natural way to
extend $\phi$ to a map, which we also call $\phi$, on the vertices of
$\B(\son)$.  If $\al$ is an arc in $\son$, we call the set of
nontrivial boundary curves of a regular neighborhood of $\al$ the
\emph{pushoffs} of $\al$.

If $\al$ is a nonseparating arc starting and ending at the same
puncture, or a separating arc with two pushoffs, then it follows from Lemmas~\ref{graphspreserved},
 \ref{linearorcyclic}, and~\ref{Npreserved} that the images under $\phi$ of the
pushoffs of $\al$ bound a punctured annulus.  Define $\phi(\al)$ to be
the unique arc in this annulus (note $\phi(\al)$ must be
topologically equivalent to $\al$).

If $\al$ is an arc with a single pushoff, then this pushoff must
be a 2-curve.  In this case, $\phi(\al)$ is determined by the image of
its pushoff and the following natural action of $\phi$ on the
punctures of $\son$: send the puncture ``shared'' by a pair of dual
2-curves to the puncture shared by the images under $\phi$ of those
2-curves (we are using Lemma~\ref{dualitypreserved}).

In his thesis, Korkmaz gives a proof that this action on the punctures of $\son$ is
well-defined \cite{mkthesis}.  He further proves:

\begin{prop}\label{korksimp} $\phi$ is a simplicial map of
$\B(\son)$.\end{prop}

Both arguments hold in our setting, since they only use properties of
superinjective maps which we have already established.


It follows easily from the injectivity of the action of $\phi$ on curves
that:

\begin{lem}\label{inj} The action of $\phi$ on arcs is
injective.\end{lem}

Now, we can finally show:

\begin{prop} $\phi$ is induced by $f$.\end{prop}

\bpf
As explained above, Proposition~\ref{korksimp} and Lemma~\ref{inj} imply that $\phi$ agrees with $f$ on all of
$\B(\son)$, and hence $\phi$ and $f$ agree on all 2-curves.  It now
follows that $\phi$ and $f$ agree on all of $\C(\son)$.
Indeed, given any curve $c$ in $\son$,
we can \emph{fill} $\son-c$ with 2-curves in the sense that $c$ is
the unique curve disjoint from the collection of 2-curves.  Since
$\phi$ and $f$ both preserve disjointness, it follows that $f$ and
$\phi$ agree on $c$.
\epf

This completes the proof of Theorem~\ref{sithm}.

\p{Remarks.} Our techniques give a particularly streamlined proof of the theorem of Bell--Margalit that any superinjective map of $\C(S_{0,n})$ ($n \geq 5$) is induced by $\mcg(S_{0,n})$.  For a genus 0 surface, a ``linear'' pants decomposition consists of a chain of punctured annuli with a twice-punctured disk on each end.  It follows then from Lemma~\ref{graphspreserved} that a superinjective map preserves the topological types of curves.  In their proof, this is the main difficult step.

We also note that the corresponding theorem of Irmak for $S_{2,0}$ follows immediately since, by the same argument of Luo used in Section~\ref{sotsec} to show
$\C(S_{1,2}) \cong \C(S_{0,5})$, we have that
$\C(S_{2,0}) \cong \C(S_{0,6})$.  It then follows from the result of Bell--Margalit that any superinjective map of $\C(S_{2,0})$ is an automorphism, and Ivanov proved that any automorphism of $\C(S_{2,0})$ is induced by an element of $\mcg(S_{2,0})$.


\bibliographystyle{plain}
\bibliography{g1}

\end{document}